\documentclass[12pt,oneside]{amsart}
\usepackage[ascii]{inputenc}
\usepackage[english]{babel}
\usepackage{amssymb}
\newcommand{\rank}{\operatorname{rank}}
\newcommand{\Lin}{\operatorname{Lin}}

\renewcommand{\Im}{\operatorname{Im}}
\newtheorem*{assum}{Assumption}
\newtheorem*{thm}{Theorem}
\newtheorem*{cor}{Corollary}
\newtheorem{Thm}{Theorem}
\newtheorem{Lem}{Lemma}
\theoremstyle{definition}
\newtheorem*{defins}{Definitions}
\newtheorem*{defin}{Definition}
\newtheorem*{rem}{Remark}

\title{On the Invariant Subspace Problem for Dissipative Operators in Krein spaces}
\author{A.~A.~Shkalikov${}^{1,2}$}

\thanks{${}^{1}$ This report
was presented as a 50-minutes lecture at Symposium on Operator
Theory in Durham (UK) held at 4-12 August, 2005 (see the
presentation at
http://maths.dur.ac.uk/events/Meetings/LMS/2005/OTSA/talks.html).
Earlier these results  were announced at the Conference on Krein
Space Operator Theory held in Berlin, December, 2004. The author
thanks the organizers for the invitation and support.}

\thanks{${}^{2}$ The work was supported by the Russian Foundation of Basic Research
(project No 04-01-00412) and by a grant of the President of
Russian Federation for the support of leading scientific schools
(project No NSc-1927.2003.1)}
\address{Department of Mechanics and Mathematics\\
Moscow State University\\
119992 Moscow, Russia} \email{shkalikov@yahoo.com}

\begin{document}
\maketitle 
\begin{abstract} We relax assumptions for a dissipative operator in Krein
space to possess a maximal non-negative  invariant subspace. Our
main result is a generalization of a well-known
Pontrjagin-Krein-Langer-Azizov theorem. Then we investigate the
semigroup properties of the restriction of the operator onto the
invariant subspace.
\end{abstract}

 \vskip 15 pt
 Let \(H\) be a separable Hilbert space and
\(J=P_+-P_-\) be a canonical symmetry (\(J^2=P_++P_-=1\)). The
space $H$ equipped with indefinite inner product
\[
    [x,y]=(Jx,y),\qquad x,y\in H
\]
is called Krein space and denoted by \(\mathcal K=\{H,J\}\) (or
Pontrjagin space \(\Pi_{\varkappa}=\{H,J\}\) if \(\rank
P_+=\varkappa<\infty\)).

\begin{defins}

1. A subspace \(\mathcal L\) is \emph{nonnegative} in \(\mathcal
K\) if \([x,x]\geqslant 0\) \(\forall x\in\mathcal L\). It is
maximal nonnegative if there are no proper extensions of
\(\mathcal L\).


2. An operator \(A\) is \emph{dissipative} in \(H\) if
\[
    \Im (Ax,x)\geqslant 0\qquad\forall x\in\mathcal D(A).
\]
It is maximal dissipative if there are no proper dissipative
extensions of \(A\) (\(\Leftrightarrow\mathbb C^-\subset\rho(A)\),
where \(\mathbb C^-\) is the open lower-half plane).

3. \(A\) is dissipative in Krein space \(\mathcal K=\{H,J\}\) if
\(JA\) is dissipative in \(H\). \(A\) is \(m\)-dissipative in
\(\mathcal K\) if \(JA\) is \(m\)-dissipative in \(H\).
\end{defins}

Symmetric and self-adjoint operators in \(\mathcal K\) are defined
analogously.
\vskip 5 pt

 Let \(H=H_+\oplus H_-\),
\(H_{\pm}=P_{\pm}(H)\), \(\mathcal D_{\pm}=\mathcal D(A) \cap
H_{\pm}\).

\begin{assum}
\(\mathcal D(A)=\mathcal D_+\oplus\mathcal D_-\) (actually, it is
sufficient to assume that \(\mathcal D_+\oplus\mathcal D_-\) is a
core of \(A\)) \(\Leftrightarrow\) \(A\) admits matrix
representation
\[
    A=\begin{pmatrix}A_{11}&A_{12}\\A_{21}&A_{22}\end{pmatrix}=
    \begin{pmatrix}P_+AP_+&P_+AP_-\\P_-AP_+&P_-AP_-\end{pmatrix},
\]
where \(x=x_++x_-\) are identified with the columns
\(x=\begin{pmatrix}x_+\\x_-
\end{pmatrix}\).
\end{assum}

\medskip
{\large\textbf{Background}}

\medskip
\begin{thm}[Sobolev, 1941, 1962]
A selfadjoint operator in \(\Pi_1\) has at least one eigenvector
corresponding to an eigenvalue \(\lambda\in\overline{\mathbb
C}^+\).
\end{thm}

\begin{thm}[Pontrjagin, 1944]
Let \(A\) be selfadjoint in \(\Pi_{\varkappa}\), \(\varkappa<\infty\). Then
\begin{itemize}
\item[(a)] \(\exists\) maximal nonnegative subspace \(\mathcal L^+\) invariant
with respect to \(A\);
\item[(b)] among these subspaces \(\exists\mathcal L^+\) such that \(\sigma
(A^+)\subset\overline{\mathbb C}^+\), \(A^+=A|_{\mathcal L^+}\).
\end{itemize}
\end{thm}

\begin{thm}[Langer, 1961]
Let \(A\) be selfadjoint in \(\mathcal K\) and
\begin{itemize}
\item[(i)] \(\mathcal D(A)\supset H_+\) (\(\Leftrightarrow A_{11}\text{ and }A_{21}\)
are bounded);
\item[(ii)] \(A_{12}\) is compact.
\end{itemize}
Then the properties (a) and~(b) hold.
\end{thm}

\begin{thm}[Krein, 1948, 1964]
Analogues of Pontrjagin and Langer theorems are true for unitary operators in
\(\Pi_{\varkappa}\) and \(\mathcal K\), respectively.
\end{thm}

M.~Krein proposed a shorther elegant approach to prove the
property (a) by means of the Shauder--Tikhonov fixed point
theorem.

\begin{thm}[Krein and Langer, 1971; Azizov, 1972]
Let \(A\) be \(m\)-dissipative in \(\Pi_{\varkappa}\). Then (a) and~(b)
hold.
\end{thm}

\begin{thm}[Azizov, Khoroshavin, 1981]
Let \(A\) be a contraction in Krein space and \(A_{12}\) be compact. Then
(a) and (b) hold if \(\mathbb C^-\) is replaced by the open unit disk.
\end{thm}

\begin{thm}[Azizov, 1985]
The analogue of the previous result holds for \(m\)-dissipative
operators in \(\mathcal K\) provided that \(\mathcal D(A)\supset
H_+\) and \(A_{12}\) is \(A_{22}\)-compact.
\end{thm}

\begin{thm}[Shkalikov, 2004]
Let
\begin{itemize}
\item[(i)] \(A\) be dissipative in \(\mathcal K\);
\item[(ii)] \(A_{22}\) be \(m\)-dissipative in \(H_-\) (\(\Leftrightarrow
\exists (A_{22}-\mu)^{-1}\) for some \(\mu\in\mathbb C^{-}\));
\item[(iii)] \(F(\mu):=(A_{22}-\mu)^{-1}A_{21}\) be bounded;
\item[(iv)] \(G(\mu):=A_{12}(A_{22}-\mu)^{-1}\) be compact;
\item[(v)] \(S(\mu):=A_{11}-A_{12}(A_{22}-\mu)^{-1}A_{21}\) be bounded.
\end{itemize}
Then (a) and~(b) hold.
\end{thm}

The  main novelty of this theorem is that the Langer condition
\(\mathcal D(A)\supset H_+\) is dropped out. In particular, for a
model matrix operator
\[
    A=\begin{pmatrix}u(x)&\dfrac{d}{dx}\\[3mm] \dfrac{d}{dx}&
    \dfrac{d^2}{dx^2}\end{pmatrix},\qquad
    J=\begin{pmatrix}1&0\\0&-1\end{pmatrix}
\]
which is selfadjoint in \(\mathcal K=\{H,J\}\), \(H=L_2[0,1]\times
L_2[0,1]\), provided that the domain of \(A\) is chosen properly,
one can guarantee the validity of properties (a) and~(b).

\textbf{The main goal of this report} is to prove the properties
(a) and (b) provided that only assumptions~(i)-(iv) are valid.

{\em It turns out that we need no assumptions for the transfer
function \(S(\mu)\)}.

New problems arise if we start working with unbounded entries and
reject Langer condition \(\mathcal D(A)\supset H_+\). In this case,
if we succeed to prove~(a) and~(b), we come to the following interesting
problems
\begin{itemize}
\item[(c)] does the operator \(A^+=A|_{\mathcal L^+}\) generate a
\(C_0\)-semigroup, or holomorphic semigroup?
\end{itemize}

We shall provide some sufficient conditions for positive answer to
this question. \vskip 8 pt

 A subspace \(\mathcal L\) is
\(A\)-invariant in classical sense if \(\mathcal L\subset\mathcal
D(A)\) and \(A(\mathcal L)\subset \mathcal L\). We accept the
following
\begin{defin}
\(\mathcal L\) is \(A\)-invariant if \(\mathcal D(A)\cap\mathcal
L\) is dense in \(\mathcal L\) and \(Ax\in\mathcal L\) for all
\(x\in \mathcal D(A)\cap\mathcal L\).
\end{defin}

Let us formulate the main results.

\begin{Thm}\label{A}
Conditions (i)--(iv) imply property (a), i.e. there exists an
$A$-invariant maximal non-negative  subspace \(\mathcal L^+\). The
property \(\mathcal L^+\subset\mathcal D(A)\) (i.e. \(\mathcal
L^+\) is $A^+$-invariant in classical sense) holds if and only if
$S(\mu)$ is bounded.
\end{Thm}

\begin{Thm}\label{B}
Property (b) holds if and only if assumption (i) is replaced by
\begin{itemize}
\item[(i')] \(A\) is \(m\)-dissipative in \(\mathcal K\).
\end{itemize}
\end{Thm}

For short we accept the following notation.
\begin{defin}
We say that the operator \(B\) is a generator of \(H_0\)-semigroup
if \(\forall\varepsilon>0\) the operator  \(B-\varepsilon\)
generates a holomorphic semigroup.
\end{defin}

Remind that $A^+$ is the restriction of $A$ onto the invariant
subspace \(\mathcal L^+\).
\begin{Thm}\label{C}
\(iA^+\) generates a \(C_0\)-semigroup of exponential type \(0\)
is one of the following conditions holds
\begin{enumerate}
\item \(A_{12}\) is compact;
\item \(-iA_{22}\) generates an \(H_0\)-semigroup.
\end{enumerate}
\end{Thm}

\begin{Thm}\label{D}
\(iA^+\) generates an exponentially stable semigroup if either
assumption~(1) or~(2) of Theorem C  holds and \(A\) is uniformly
dissipative in \(\mathcal K\).
\end{Thm}

\begin{Thm}\label{E}
There is \(\mu\in\mathbb C^+\) such that \(iS(\mu)\) generates an
\(H_0\)-semigroup. Then \(iA^+\) generates \(H_0\)-semigroup.
\end{Thm}

\medskip
{\large\textbf{Proof of the first two theorems}}

\medskip
Assumptions~(ii)--(iv) allow to use Frobenious-Shur factorization
\[
    A-\mu=\begin{pmatrix}1&G\\0&1\end{pmatrix}
    \begin{pmatrix}S-\mu&0\\0&A_{22}-\mu\end{pmatrix}
    \begin{pmatrix}1&0\\F&1\end{pmatrix}
\]
where \(G=G(\mu)\), \(F=F(\mu)\) and \(S=S(\mu)\) is the transfer function
defined on the domain \(\mathcal D(S)=\mathcal D_+\).

\begin{Lem}\label{L1}
\[
    JA+\mu=J\begin{pmatrix}1&G\\0&1\end{pmatrix}
    \begin{pmatrix}S+\mu&0\\0&A_{22}-\mu\end{pmatrix}
    \begin{pmatrix}1&0\\F&1\end{pmatrix}
\]
\end{Lem}

\textit{Proof} is obtained by direct verification.

\begin{Lem}\label{L2}
\(\forall\mu\in\mathbb C^+\) and \(\forall x\in\mathcal D_+\) we have
\[
    (Sx_+,x_+)=\left(JA\begin{pmatrix}x_+\\-Fx_+\end{pmatrix},
    \begin{pmatrix}x_+\\-Fx_+\end{pmatrix}\right)+\mu(Fx_+,Fx_+).
\]
\end{Lem}

\textit{Proof} is obtained by direct verification.

\begin{cor}[important]
\(S=S(\mu)\) with domain \(\mathcal D(S)=\mathcal D_+\) is dissipative in
\(H_+\) provided that assumption~(i) holds. Also, \(S\) is closable. The
closure of \(S\) is \(m\)-dissipative in \(H_+\) \(\Leftrightarrow\)
\(A\) is \(m\)-dissipative in \(\mathcal K\).
\end{cor}

\begin{Lem}[important]\label{L3}
Let a subspace \(\mathcal L\) have a representation of the form
\[
    \mathcal L=\{x\;:\;x=\begin{pmatrix}x_+\\Kx_+\end{pmatrix},\;
    x_+\in H_+\}
\]
where \(K:H_+\to H_-\) is a bounded operator. Then \(\mathcal L\)
is \(A\)-invariant \(\Leftrightarrow\)
\[
    (1-KG)(A_{22}-\mu)(F+K)=K(S-\mu)
\]
(the so-called Riccati equation for \(K\)).
\end{Lem}

\begin{proof}
For \(x_+\in\mathcal D_+\)
\[
    (A-\mu)\begin{pmatrix}x_+\\Kx_+\end{pmatrix}=\begin{pmatrix}
    (S-\mu)x_++G(A_{22}-\mu)(F+K)x_+\\ (A_{22}-\mu)(F+K)x_+
    \end{pmatrix}.
\]
Assuming that \(\mathcal L\) is \(A\)-invariant we find \(y_+\in H_+\)
such that
\begin{gather*}
    [(S-\mu)+G(A_{22}-\mu)(F+K)]x_+=y_+,\\
    (A_{22}-\mu)(F+K)x_+=Ky_+.
\end{gather*}
Substituting the first equality in the second one we come to
Riccati equation for \(K\).

Conversely, Riccati equation for \(K\) implies the last two equations
with some \(y_+\), therefore the graph subspace \(\mathcal L\) is
\(A\)-invariant.
\end{proof}

\begin{rem}
Pontrjagin used the following version of Lemma 3: \(\mathcal L\)
is \(A\)-invariant \(\Leftrightarrow\)
\[
    A_{21}+A_{22}K-KA_{11}-KA_{12}K=0.
\]
However this form of Riccati equation is inconvenient while
working with unbounded entries \(A_{ij}\).
\end{rem}

\begin{Lem}\label{L4}
Assume that \(G(\mu)\) is compact for some \(\mu\in\mathbb C^+\). Then it is
compact for all \(\mu\in\mathbb C^+\) and \(\|G(\mu)\|\to 0\) as \(\mu\to
\infty\) and \(\mu\in\Lambda_{\varepsilon}^+\).
\end{Lem}

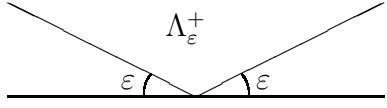
\begin{figure}[h]
\unitlength=0.1cm
\begin{picture}(50,20)
\put(0,0){\line(1,0){50}}
\put(25,0){\line(-2,1){25}}
\put(25,0){\line(2,1){25}}
\put(21,8){\text{\(\Lambda_{\varepsilon}^+\)}}
\put(15,1){\text{\(\varepsilon\)}}
\put(33,1){\text{\(\varepsilon\)}}
\qbezier(18,0)(18,2)(19,3)
\qbezier(32,0)(32,2)(31,3)
\end{picture}
\end{figure}

\textit{Proof} is simple.

\begin{Lem}\label{L5}
A subspace \(\mathcal L\) is maximal nonnegative
\(\Leftrightarrow\) \(\mathcal L\) has the graph representation
\[
    \mathcal L=\{x=\begin{pmatrix}x_+\\ Kx_+\end{pmatrix},\quad
    x_+\in H_+\}
\]
with the angle operator \(K\), \(\|K\|\leqslant 1\).
\end{Lem}

\begin{cor}
Take \(\mu\in\mathbb C^+\) such that \(\|G(\mu)\|<1/2\). Then (a)
holds  \(\Leftrightarrow\) \(\exists\) there is a contraction
\(K\) such that
\[
    F+K=(A_{22}-\mu)^{-1}(1-KG)^{-1}K(S-\mu).
\]
\end{cor}

\begin{Lem}\label{L6}
Denote \(H_S=\mathcal D(\overline{S})\subset H_+\) where
\(\overline{S}\) is the closure of \(S\) and the norm in \(H_S\) is
defined by
\[
    \|x_+\|_{H_S}=\sqrt{\|\overline{S}x_+\|^2+\|x_+\|^2}.
\]
Then there is a complete orthogonal system
\(\{\varphi_k\}_1^{\infty}\) in \(H_+\) such that
\(\{\varphi_k\}_1^{\infty}\) is a Riesz basis in \(H_S\).
\end{Lem}

\begin{proof}
If \(H_S\) is compactly embedded in \(H_+\) we take
\(\{\varphi_k\}_1^{\infty}\) consisting of eigenvectors of \(S^*
\overline{S}\). In general case additional work is required,
however, this work is routine.
\end{proof}

\begin{proof}[Proof of Theorem~A]
Let \(P_n\) be orthogonal projectors onto \(\Lin\{\varphi_k\}_1^n\) in
\(H_+\). Then \(P_n\to 1\) in \(H_+\) and \(P_n\to 1\) in \(H_S\).

Consider
\[
    A_n=\begin{pmatrix}P_nA_{11}P_n&P_nA_{12}\\ A_{21}P_n&A_{22}
    \end{pmatrix}\qquad\text{in } H_n^+\oplus H^-,
    \; H_n^+=P_n(H^+).
\]
Then \(A_n\) is \(m\)-dissipative in Pontrjagin space \(\Pi_n\) and
due to Krein--Langer--Azizov theorem~(a) holds. This implies (Lemma~3)
that
\begin{equation}\label{E1}
    F_n+K_n=(A_{22}-\mu)^{-1}(1-K_nG)^{-1}K_n(S_n-\mu).
\end{equation}
It is known that the unit ball of a separable Hilbert space is
weakly compact. The sequence of the operators $K_n$ is bounded,
hence, we can choose a weakly convergent subsequence
\(K_{n_j}\rightharpoonup K\). For short we omit the index $j$.
Since \(\|K_n\|\leqslant 1\), we have \(\|K\|\leqslant 1\). Then
\[
    F_n=FP_n\to F, \quad K_nG\Rightarrow KG \ \text{and}\ (1-K_nG)^{-1}\Rightarrow 1-KG
\]

(we essentially use here that \(G\) is compact!).

Further,
\begin{gather*}
    K_nS_n=K_nSP_n,\\
    \overline{S}P_nx\to\overline{S}x\qquad\forall x\in
    \mathcal D(\overline{S}),
\end{gather*}
Hence, \(K_nSP_nx\rightharpoonup KSx\). Therefore we can pass to the
weak limit in the equation~\eqref{E1} and obtain
\[
    F+K=(A_{22}-\mu)^{-1}(1-KG)^{-1}K(S-\mu)
\]
and by virtue of Lemma~3 property (a) holds.
\end{proof}

Let \(A^+=\overline{A}|_{\mathcal L^+}\). How to prove the
property
\begin{itemize}
\item[(b)] \(\exists\mathcal L^+\) such that \(\sigma(A^+)\subset
\overline{\mathbb C^+}\)?
\end{itemize}

We have
\[
    (\overline{A}-\mu)\begin{pmatrix}x_+\\ Kx_+\end{pmatrix}=
    \begin{pmatrix}(\overline{S}-\mu+GL)x_+\\ Lx_+\end{pmatrix},
\]
where \(L:=(A_{22}-\mu)(F+K)\), \(\mathcal D(L)=
\mathcal D(\overline S)\).

Consider
\[
    Q:\mathcal L\to H_+\quad\text{defined by } Q\begin{pmatrix}
    x_+\\ Kx_+\end{pmatrix}=x_+.
\]
\(Q\) is bounded and boundedly invertible, \(\|Q^{-1}\|\leqslant 2\).

We have
\[
    \overline{A}|_{\mathcal L^+}=Q^{-1}(S+GL)Q=Q^{-1}[1+G(1-KG)^{-1}
    K(\overline{S}-\mu)]Q,
\]
hence
\begin{gather}\label{E2}
    (\overline{A}-\alpha)|_{\mathcal L^+}=Q^{-1}[1+T(\alpha)]
    (\overline{S}(\mu)-\alpha)Q,\\ \intertext{where}\notag
    T(\alpha)=G(1-KG)^{-1}K(\overline{S}-\mu)(\overline{S}-\alpha)^{-1}
\end{gather}
is a holomorphic operator function whose values are compact operators.
Here we assumed that \((\overline{S}-\alpha)^{-1}\) exists
\(\Leftrightarrow\) \(\overline{S}\) is \(m\)-dissipative in \(H_+\)
\(\Leftrightarrow\) \(\overline{A}\) is \(m\)-dissipative in \(\mathcal K\).
It can be shown that \(\|T(\alpha)\|\to 0\) as \(\alpha\to\infty\) along
negative imaginary axis, therefore \(1+T(\alpha)\) has only discrete
spectrum in \(\mathbb C^-\).

Now we shall use the following
\begin{Lem}
\(\Im[Ax_0,x_0]=(\Im\alpha_0)\, [x_0,x_0]\) if
\(Ax_0=\alpha_0x_0\).
\end{Lem}

Therefore, all eigenvectors of \(A\) corresponding to the
eigenvalues from \(\mathbb C^-\) are of negative type, provided
that \(A\) is strictly dissipative in \(\mathcal K\). Hence, these
vectors do not belong to \({\mathcal L}^+\). This means that the
eigenvalues from \(\mathbb C^-\) do not belong to the spectrum of
$A^+$, therefore \(\sigma(A^+)\in \overline{\mathbb C}^+\).

This proves Theorem~B if we assume in addition that \(A\) (or
$A^+$) is strictly dissipative in \(\mathcal K\).

If not, we consider
\[
    A_{\varepsilon}=A+i\varepsilon P_+,\qquad\varepsilon>0.
\]
Assertion~(a) is valid for \(A_{\varepsilon}\), and it does not
have spectrum in \(\mathbb C^-\), since
\[
    \Im\left[(A+i\varepsilon P_+)\begin{pmatrix}x_+\\ Kx_+\end{pmatrix},
    \begin{pmatrix}x_+\\ Kx_+\end{pmatrix}\right]\geqslant
    \varepsilon(x_+,x_+).
\]
Write Riccati equation for \(A_{\varepsilon}\):
\[
    F+K_{\varepsilon}=(A_{22}-\mu)^{-1}(1-K_{\varepsilon}G)^{-1}
    K_{\varepsilon}(S+i\varepsilon-\mu).
\]
Take \(\varepsilon_n\to 0\) and \(K_{\varepsilon_n}=:K_n\rightharpoonup
K\).

We have
\begin{gather*}
    A_{\varepsilon}^+=Q^{-1}[1+T_{\varepsilon}(\alpha)](S+i\varepsilon-
    \alpha)Q\\ \intertext{and}
    T_{\varepsilon}(\alpha)=G(1-K_{\varepsilon}G)^{-1}K_{\varepsilon}
    (S+i\varepsilon-\mu)(S+i\varepsilon-\alpha)^{-1}\Rightarrow
    T(\alpha).
\end{gather*}
Since \(1+T_{\varepsilon}(\alpha)\) is a holomorphic operator
function of Fredholm type in \(\mathbb C_-\), and boundedly
invertible \(\forall\alpha\in \mathbb C^-\), so is
\(1+T(\alpha)\).\hfill\(\Box\) (It is important here that
$T_{\varepsilon}(\alpha)$ converges to $T(\alpha)$
 uniformly, and $1+T(\alpha)$ has discrete spectrum in $\mathbb
 C^-$).

Theorems~C--E are proved by analyzing representation~\eqref{E2}.
\vskip 7 pt

 For convenience we present here references related to
the background of the problem.

\end{document}